\documentclass{elsart} \usepackage{amsfonts}
\usepackage{amsmath}
\usepackage[all]{xy}
\usepackage{epsf}
\usepackage{epsfig}

\newcommand{\C}{{\mathbb C}}

\newcommand{\Z}{{\mathbb Z}}
\newcommand{\N}{{\mathbb N}}
\newcommand{\bE}{{\mathbb E}}
\newcommand{\bB}{{\mathbb B}}

\newcommand{\cB}{{\mathcal B}}
\newcommand{\cC}{{\mathcal C}}

\newcommand{\isom}{\cong}
\newcommand{\D}[1]{\ensuremath{D^{(#1)}}}
\newcommand{\T}[1]{\ensuremath{\mathrm{T}(#1)}}

\newcommand{\set}[2]{\ensuremath{\{#1|#2\}}}
\newcommand{\Hom}{\ensuremath{\mathrm{Hom}}}

\newcommand{\Multi}{\ensuremath{\mathit{Multi}}}
\newcommand{\Sing}{\ensuremath{\mathrm{Sing}}}
\newcommand{\Multibin}{\ensuremath{\mathit{Multi}}_{\epsfig{file=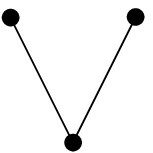,height=2mm}}}

\newcommand{\h}{\mbox{-}}

\newcommand{\tens}{\mbox{$\otimes$}}

\newcommand{\Rmod}{\ensuremath{\mbox{$R\mbox{-}\mathit{Mod}$}}}
\newcommand{\Hmod}{\ensuremath{\mbox{$H\mbox{-}\mathit{Mod}$}}}

\renewcommand{\H}[1]{\ensuremath{H^{\otimes #1}}}

\newcommand{\bottom}{\bot}

\newcommand{\acts}{\mbox{$\cdot$}}

\newcommand{\ps}[1]{[\![#1]\!]}

\newcommand{\bemph}[1]{\emph{\textbf{#1}}}
\newcommand{\definition}[1]{\textbf{#1}}

\newcommand{\revar}[2]{\xymatrix{*=0{} \ar@<2pt>[r]^{#1} & *=0{}
\ar@<2pt>[l]^{#2}}} 
\newcommand{\revmapsto}[2]{\xymatrix{*=0{} \ar@{|->}@<3pt>[r]^{#1} & *=0{}
\ar@{|->}@<3pt>[l]^{#2}}} 

\newcommand{\genericblank}{\vcenter{
\xymatrix{ & & \ldots &  \\
	   & & *=0{} \ar@{-}[llu] \ar@{-}[lu] \ar@{-}[ru]}}}

\newcommand{\generic}[3]{\vcenter{
\xymatrix{{#1}_1&{#1}_2 &  \ldots	& {#1}_{#3}  \\
		& 	& *=0{} \ar@{-}[llu]^{{#2}_1}
			  \ar@{-}[lu]_{{#2}_2}\ar@{-}[ru]_{{#2}_{#3}}}}}

\newcommand{\genericnolabel}[3]{\vcenter{
\xymatrix{{#1}_1&{#1}_2 &  \ldots	& {#1}_{#2}  \\
		& 	& *=0{\txt{\\ \\{$#3$}}} \ar@{-}[llu] 
			  \ar@{-}[lu]\ar@{-}[ru]}}}

\newcommand{\flatmodule}[4]{\vcenter{
\xymatrix{{#1}_1&{#1}_2 &  \ldots	& {#1}_{#3}  \\
		& 	& *=0{\txt{\\ \\{$#4$}}} \ar@{-}[llu]^{{#2}_1}
			  \ar@{-}[lu]_{{#2}_2}\ar@{-}[ru]_{{#2}_{#3}}}}}

\newcommand{\flattwoleafed}[5]{\vcenter{
\xymatrix@C=10pt@R=15pt{{#1} &  & {#2}  \\
	   & *=0{\txt{\\ \\{$#3$}}} \ar@{-}[lu]^{#4} \ar@{-}[ru]_{#5}}}}

\newcommand{\flatthreeleafed}[4]{\vcenter{
\xymatrix@C=10pt@R=15pt{{#1} &{#2} & {#3} \\
	   & *=0{\txt{\\ \\{$#4$}}} \ar@{-}[lu] \ar@{-}[u] \ar@{-}[ru]}}}

\newcommand{\flatthreeleafedlabelled}[7]{\vcenter{
\xymatrix@C=10pt@R=15pt{{#1} &{#2} & {#3} \\
	   & *=0{\txt{\\ \\{$#4$}}} \ar@{-}[lu]^{#5} \ar@{-}[u]_{#6} \ar@{-}[ru]_{#7}}}}

\newcommand{\threeleafedleft}[8]{\vcenter{
\xymatrix@C=10pt@R=15pt{ {#1}&                               &{#2} & &  \\
	       & *=0{} \ar@{-}[lu]^{#5} \ar@{-}[ru]_{#6} &     & {#3}&  \\
               & & *=0{\txt{\\ \\{$#4$}}} \ar@{-}[lu]^{#7} \ar@{-}[ru]_{#8}}}}

\newcommand{\threeleafedright}[8]{\vcenter{
\xymatrix@C=10pt@R=15pt{     &                               &{#2} & &{#3}&  \\
	   &{#1}& & *=0{} \ar@{-}[lu]^{#5} \ar@{-}[ru]_{#6} &  \\
               & & *=0{\txt{\\ \\{$#4$}}} \ar@{-}[lu]^{#7} \ar@{-}[ru]_{#8}}}}

\newcommand{\fourleafdoubleA}{\vcenter{
\xymatrix@C=10pt@R=20pt{ V &  & V & & V & & V  \\
	     & *=0{} \ar@{-}[lu]^{x} \ar@{-}[ru] 
                &   & &   & *=0{} \ar@{-}[lu]^{y} \ar@{-}[ru] &  \\
             &  &   & *=0{\txt{\\ \\$V$}} \ar@{-}[llu]^{z} \ar@{-}[rru]}}}

\newcommand{\fourleafdoubleB}[5]{\vcenter{
\xymatrix@C=10pt@R=20pt{ {#1} &  & {#2} & & {#3} & & {#4}  \\
	     & *=0{} \ar@{-}[lu] \ar@{-}[ru] 
                &   & &   & *=0{} \ar@{-}[lu] \ar@{-}[ru] &  \\
             &  &   & *=0{\txt{\\ \\{$#5$}}} \ar@{-}[llu] \ar@{-}[rru]}}}

\newcommand{\fourleafdoubleC}[5]{\vcenter{
\xymatrix@C=10pt@R=20pt{ {#1}_1 &  & {#1}_2 & & {#1}_3 & & {#1}_4  \\
	     & *=0{} \ar@{-}[lu]^{{#3}_1} \ar@{-}[ru]_{{#3}_2}
                &   & &   & *=0{} \ar@{-}[lu]^{{#4}_1} \ar@{-}[ru]_{{#4}_2} &  \\
             &  &   & *=0{\txt{\\ \\{$#2$}}} \ar@{-}[llu]^{{#5}_1} \ar@{-}[rru]_{{#5}_2}}}}

\newcommand{\twoleavedrighttwo}[6]{\vcenter{
\xymatrix@C=10pt@R=15pt{    & & {#2}& &  \\
	                {#1}& & *=0{} \ar@{-}[u]^{#6} &  \\
                            & *=0{\txt{\\ \\{$#3$}}} \ar@{-}[lu]^{#4} \ar@{-}[ru]_{#5}}}}

\newcommand{\twoleafedtail}[6]{\vcenter{
\xymatrix@C=10pt@R=15pt{ {#1}&                                         &{#2}  \\
	                     & *=0{} \ar@{-}[lu]^{#4} \ar@{-}[ru]_{#5} &       \\
                             & *=0{\txt{\\ \\{$#3$}}} \ar@{-}[u]^{#6} }}}

\newcommand{\flatthreeleafedleft}[4]{\vcenter{
\xymatrix@C=10pt@R=15pt{ {#1}&                               &{#2} &{#3}&  \\
	       & *=0{} \ar@{-}[lu] \ar@{-}[ru] &     & *=0{} \ar@{-}[u]&  \\
               & & *=0{\txt{\\ \\$#4$}} \ar@{-}[lu] \ar@{-}[ru]}}}

\newcommand{\flatthreeleafeddown}{\vcenter{
\xymatrix@C=10pt@R=15pt{ & &  \\
	   & *=0{} \ar@{-}[lu] \ar@{-}[u] \ar@{-}[ru] \\
           & *=0{} \ar@{-}[u]  }}}

\newcommand{\genericA}[3]{\vcenter{
\xymatrix@C=10pt@R=15pt{{#1}_1&{#1}_2 &  \ldots	& {#1}_{#3}  \\
		& 	& *=0{\txt{\\ \\$#2$}} \ar@{-}[llu]
			  \ar@{-}[lu]\ar@{-}[ru]}}}

\newcommand{\generictailA}[3]{\vcenter{
\xymatrix@C=10pt@R=15pt{
{#1}_1	&{#1}_2	&\ldots	& {#1}_{#3}  \\
	&	&*=0{} \ar@{-}[llu]
			  \ar@{-}[lu] \ar@{-}[ru] \\
	&	&*=0{\txt{\\ \\$#2$}} \ar@{-}[u]}}}

\newcommand{\generictopA}[3]{\vcenter{
\xymatrix@C=10pt@R=15pt{
	&	&	&{#1}_k			&	&\\
{#1}_1	&{#1}_2	&\ldots	&*=0{} \ar@{-}[u]	&\ldots	&{#1}_{#3}  \\
	&	& *=0{\txt{\\ \\$#2$}} \ar@{-}[llu]
			  \ar@{-}[lu] \ar@{-}[ru] \ar@{-}[rrru]}}}

\newcommand{\genericcompositionB}{\vcenter{
\xymatrix@C=10pt@R=15pt{
A_{1}	& A_{2}	&\ldots	& A_{m}&\\
	&	& *=0{} \ar@{-}[llu] \ar@{-}[lu]\ar@{-}[ru]
			&A_{m+1}	&\ldots & A_n\\ 
	&		&&&*=0{\txt{\\ \\$B$}} \ar@{-}[llu]\ar@{-}[lu]
			 \ar@{-}[ru]  }}}

\newcommand{\genericcompositionA}{\vcenter{
\xymatrix@C=10pt@R=15pt{&& A_{i_1}& A_{i_2}	&\ldots	& A_{i_m}&\\
	  A_1&A_2&	&\ldots	& *=0{} \ar@{-}[llu]
			  \ar@{-}[lu]\ar@{-}[ru]&\ldots &A_n\\
	  &&	&	&*=0{\txt{\\ \\$B$}} \ar@{-}[llllu]
			  \ar@{-}[lllu] \ar@{-}[u]\ar@{-}[rru]  }}}

\newcommand{\flatthreeleafedleftsides}[6]{\vcenter{
\xymatrix{&                    & & &  \\
	  & *=0{} \ar@{-}[lu]^{#1} \ar@{-}[ru]_{#2} & & *=0{} \ar@{-}[u]_{#3}&  \\
          & & *=0{\txt{\\ \\{$#6$}\\}} \ar@{-}[lu]^{#4} \ar@{-}[ru]_{#5}}}}

\newcommand{\flatthreeleafedleftsidestop}[6]{\vcenter{
\xymatrix{A_1&                    &A_2 &A_3 &  \\
	  & *=0{} \ar@{-}[lu]^{#1} \ar@{-}[ru]_{#2} & & *=0{} \ar@{-}[u]_{#3}&  \\
          & & *=0{\txt{\\ \\{$#6$}\\}} \ar@{-}[lu]^{#4} \ar@{-}[ru]_{#5}}}}

\newcommand{\genericcomposition}{\vcenter{
\xymatrix{&&	&	&\ldots	& &\\
	  &&	&\ldots	& *=0{} \ar@{-}[llu]^{{x}_1}
			  \ar@{-}[lu]_{{x}_2}\ar@{-}[ru]_{{x}_{n}}&\ldots &\\
	  &&	&	&*=0{} \ar@{-}[llllu]^{{y}_1}
			  \ar@{-}[lllu]_{{y}_2} \ar@{-}[u]_{y_k}
\ar@{-}[rru]_{{y}_{m}}  }}}

\newcommand{\generictail}[3]{\vcenter{
\xymatrix{
{#1}_1	&{#1}_2	&\ldots	& {#1}_{#3}  \\
	&	&*=0{} \ar@{-}[llu]^{{#2}_1-y}
			  \ar@{-}[lu]_{{#2}_2-y} \ar@{-}[ru]_{{#2}_{#3}-y} \\
	&	&*=0{} \ar@{-}[u]^{y}}}}

\newcommand{\oldgenerictail}[3]{\vcenter{
\xymatrix{
{#1}_1	&{#1}_2	&\ldots	&{#1}_k	&\ldots	& {#1}_{#3}  \\
	&	& 	&*=0{} \ar@{-}[lllu]^{{#2}_1-{#2}_k}
			  \ar@{-}[llu]_{{#2}_2-{#2}_k} \ar@{-}[u]_{0}
\ar@{-}[rru]_{{#2}_{#3}-{#2}_k} \\
	&	&	&*=0{} \ar@{-}[u]^{{#2}_k}}}}

\newcommand{\generictop}[3]{\vcenter{
\xymatrix{
	&	&	&{#1}_k			&	&\\
{#1}_1	&{#1}_2	&\ldots	&*=0{} \ar@{-}[u]^{y}	&\ldots	&{#1}_{#3}  \\
	&	& *=0{} \ar@{-}[llu]^{{#2}_1}
			  \ar@{-}[lu]_{{#2}_2} \ar@{-}[ru]^{{#2}_k-y}
\ar@{-}[rrru]_{{#2}_{#3}}}}}

\newcommand{\genericflattop}[3]{\vcenter{
\xymatrix{
{#1}_1	&{#1}_2	&\ldots	&{#1}_k			&\ldots	&{#1}_{#3}  \\
	&	& *=0{} \ar@{-}[llu]^{{#2}_1}
			  \ar@{-}[lu]_{{#2}_2} \ar@{-}[ru]^{{#2}_k}
\ar@{-}[rrru]_{{#2}_{#3}}}}}

\newcommand{\threenleaf}{
\xymatrix{ 
a_1	& a_2 	& a_3 	& a_4 	& a_5	&  \ldots	& a_n  \\
	&*=0{}\ar@{-}[lu]_{y_1} \ar@{-}[u]_{y_2}
		&	& *=0{}\ar@{-}[lu]_{y_3}\ar@{-}[u]_{y_4}\ar@{-}[ru]_{y_5}
				&\ldots	& *=0{}\ar@{-}[ru]_{y_n} \\
	&	&	&	&*=0{}\ar@{-}[lllu]^{x_1} \ar@{-}[lu]_{x_2}
				\ar@{-}[ru]_{x_n}  
} }

\newcommand{\threeleafedleftnoedgelabels}[5]{\vcenter{
\xymatrix@C=10pt@R=15pt{ {#1}&                               &{#2} & &  \\
	       & *=0{\txt{\\ \\{$#4$}}} \ar@{-}[lu] \ar@{-}[ru] &     & {#3}&  \\
               & & *=0{\txt{\\ \\{$#5$}}} \ar@{-}[lu] \ar@{-}[ru]}}}

\newcommand{\threeleafedrightnoedgelabels}[5]{\vcenter{
\xymatrix@C=10pt@R=15pt{     &                               &{#2} & &{#3}&  \\
	   &{#1}& & *=0{\txt{\\ \\{$#4$}}} \ar@{-}[lu] \ar@{-}[ru] &  \\
               & & *=0{\txt{\\ \\{$#5$}}} \ar@{-}[lu] \ar@{-}[ru]}}}

\begin{document}
\begin{frontmatter}
\title{Relaxed multicategory structure of a global category of rings and modules} 
\author{Craig T. Snydal\thanksref{email}}
\address{Department of Pure Mathematics and Mathematical Statistics \\
Centre for Mathematical Sciences\\
University of Cambridge\\
Wilberforce Road\\
Cambridge CB3 0WB}
\thanks[email]{ctsnydal@dpmms.cam.ac.uk}
\date{31 January 2001}
\begin{abstract}
In this paper we describe how to give a particular global category of rings
and modules the structure of a relaxed multicategory, and we describe an
algebra in this relaxed multicategory such that vertex algebras appear as
such algebras.
\end{abstract}

\end{frontmatter}

Our intention for this paper is to describe a method for giving a category of
modules for a cocommutative Hopf algebra, the structure of a relaxed
multicategory.  Relaxed multicategories are a generalization of Lambek's
multicategories \cite{lambek}, and were introduced by Richard Borcherds
in the paper, \cite{bor}, as the natural setting for vertex algebras.  The
idea was to give a category of modules for a Hopf algebra enough extra
structure that vertex algebras would arise naturally as monoids.

It is worth mentioning here that although relaxed multicategories bear a
strong resemblance to colored operads, they are nonetheless very different.
Beilinson and Drinfeld have used used multicategories/colored operads to
investigate chiral algebras \cite{beilinson_drinfeld}, and recently Soibelman
and Kontsevich looked further into this approach \cite{soibelman_kont}, but
it is fundamentally different from the relaxed multicategory treatment.

In the treatment that follows, we begin by reviewing the definition of a
relaxed multicategory.  We then define a global category of rings and modules
and show that the forgetful functor to the category of rings defines a
bifibration.  Next we describe what types of singularities we will be working
with, and we define binary singular multimaps.  We then go on to define more
general singular maps by using pushouts and pullbacks in our global category
of rings and modules.  Finally we show that suitable collections of these
maps provide the structure of a relaxed multicategory, and we define an
algebra in this setting.  These algebras give a natural interpretation of the
``locality'' axiom of traditional vertex algebras, and they formalize the
notion of operator product expansion.

\section{Definition of a relaxed multicategory}
In order to give the definition of a relaxed multicategory, we will need to
work over a certain category of trees.  The definition we take for our
category of trees is due to Tom Leinster \cite{leinster_structures},
\cite{leinster_enrichment} and seems to be a natural one arising from higher
dimensional categorical considerations.  Other categories of trees have been
defined by Ginzburg and Kapranov \cite{ginzburg_kapranov}, Manin
\cite{manin}, and Soibelman \cite{soibelman} which differ slightly from this
definition in both their collections of objects and their allowed maps.

For each natural number, $n$, (including zero) we define, $\T{n}$, the
\definition{set of $n$\h leafed trees} by the recursion:
\begin{enumerate}
\item For some formal symbol, $\bullet$, we have $\bullet \in \T{1}$ and 
\item For natural numbers $n, k_1, \ldots, k_n$, and for $t_i \in \T{k_i}$, we
have $\langle t_1, \ldots, t_n \rangle \in \T{k_1+ \cdots k_n}$.  
\end{enumerate}
We may represent these trees pictorially as, for example, $\bullet = \bullet$, 
$\langle \bullet \rangle = \epsfig{file=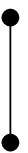,height=4mm}$, 
$\langle \bullet, \ldots,\bullet \rangle = \epsfig{file=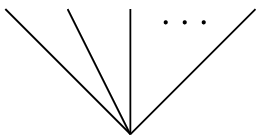,height=4mm}$, 
$\langle \langle \bullet,\bullet \rangle, \bullet \rangle = \epsfig{file=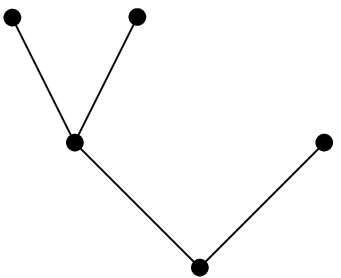,height=4mm}$.
We have $\langle \ \rangle \in \T{0}$, and we represent this empty
tree by $\circ$.  In $\T{0}$ we also have trees of the form $\langle \langle \ \rangle,
\langle \ \rangle\rangle$ which we consider to be trees with no leaves, and
which are represented as $\epsfig{file=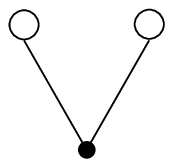,height=5mm}$

In this pictorial representation of the trees, the bottom vertex (or node) is
called the \definition{root} of the tree, and among the other vertices, those
which are joined to exactly one edge (excluding the root) are called the
\definition{external vertices} or the \definition{leaves}.  The remaining
vertices are \definition{internal vertices}.  The \definition{level of a
vertex} is defined to be the number of edges separating that vertex and the
root.  If all the leaves of a tree have level one, then the tree is called
\definition{flat} or a \definition{corolla}.

Trees compose by grafting a root to a specified leaf, and these compositions
define associative maps $\T{n} \times \T{m} \rightarrow \T{m+n-1}$ for
natural numbers $n$ and $m$ ($n >0$).  For any two trees $p, q$ with the same
number of leaves, we say that $p$ is a \definition{refinement} of $q$ if $p$
arises after a succession (possibly zero) of the following moves:
\begin{enumerate}
\item a vertex is replaced by an edge (i.e., any subtree $t$ of a given tree
can be replaced with $\langle t \rangle$),
\item any proper subtree, $\langle t \rangle$, of a given tree is replaced
with $t$ (i.e., this mostly means that an edge which is not connected to a
leaf can be shrunk down to a vertex).
\end{enumerate}
We give $\T{n}$ the structure of a category by defining a unique morphism $q
\rightarrow p$ whenever $p$ is a refinement of $q$.

\begin{note}
This is the same as the usual maps of trees, $Vert(T) \rightarrow
Vert(T^\prime)$ respecting 
\[\xymatrix{Edge(T) \ar@<.5ex>[r]^s \ar@<-.5ex>[r]_t& Vert(T)}.\]
\end{note}

We are now ready to define a relaxed multicategory.

\begin{defn}\label{defn_rel_multicat}
A \definition{relaxed multicategory} consists of a collection of objects,
$\cC$, together with a collection of \definition{multimaps} from $A_1,
\ldots, A_n$ to $B$ for any $n+1$ objects, $A_1, \ldots, A_n, B$ in $\cC$,
and any $n$\h leafed tree, $p$.  This collection is denoted 
\[\Multi_p(A_1, \ldots, A_n; B).\]
This data satisfies the following axioms:
\begin{enumerate}
\item \bemph{Identities:} For each object $A \in \cC$, there is a unique
identity map $1_A \in \Multi_\bullet(A;A)$.
\item \bemph{Composition:} Inherited from the grafting of trees, given trees 
$p \in \T{n}$ and $q \in \T{m}$ ($m\geq 1$) and objects $A_1, \ldots, A_n,
B_1,\ldots B_m, C \in \cC$, we have a map
\[\begin{split}\Multi_q(B_1, &\ldots, B_m; C) \tens \Multi_p(A_1, \ldots, A_n; B_i) \\
&\longrightarrow \Multi_{q\circ_i p}(B_1, \ldots, A_1, \ldots, A_n, \ldots, B_m; C),
\end{split}\]
where $q\circ_i p$ is the tree arrived at by grafting the root of the tree
$p$ to the $i$th leaf of the tree $q$.  This composition is associative and
must agree with identities on objects.
\item \bemph{Refinement:} Maps between trees, $p \rightarrow q$, induce maps between
multimaps in the opposite direction,
\[\Multi_q \longrightarrow \Multi_p \]
which are natural with respect to composition.
\end{enumerate}
A relaxed multicategory has an underlying category whose morphisms are given
by $\Hom(A,B) = \Multi_\bullet(A;B)$.
\end{defn}

\section{The Global Category of Rings and Modules}
The category which will provide the natural setting for working with vertex
algebras will be a the global category of rings and modules.  The intuitive
idea is that we want to use the machinery of limits and colimits to give a
certain category of modules some extra multicategory structure.  We will need
to complicate the situation slightly by giving our rings and modules the
structure of modules for a cocommutative Hopf algebra.

Fix a commutative ring $R$ with unit and take $\Rmod$ to be the symmetric
monoidal closed category of $R$\h modules.  We then take $H$ to be a
cocommutative Hopf algebra object in $R$\h modules.  Recall that a Hopf
algebra is a module, $H$, over a commutative ring, $R$ (with unit), that has
both the structure of an algebra and a coalgebra where the algebra and
coalgebra maps are compatible with one another (i.e. the maps giving $H$ the
structure of an algebra are maps of coalgebras, and vice versa).  A Hopf
algebra also possesses a multiplication and comultiplication reversing
bialgebra map, $S:H^{op} \rightarrow H$, called antipode (see \cite{abe}).

Since $H$ is a monoid, we can form the category, $\Hmod = [H, \Rmod]$, of $H$\h
modules.  This category has tensor products inherited from $R$\h modules, and
the cocommutativity and coassociativity of $H$ give $\Hmod$ the structure of
an enriched symmetric monoidal category (with unit $R$).  It can be easily
checked that the closed structure of the category of $R$\h modules carries
over to $\Hmod$.  Because the category of $R$\h modules is complete and
cocomplete, it follows from some basic results of enriched category theory
(see \cite{kelly}) that the enriched presheaf category $\Hmod$ is complete
and cocomplete, and the limits and colimits are computed pointwise.

We now want to consider the global category, $\bE$, of rings and modules
built up from $\Hmod$.  It has as objects, pairs $(L,A)$ where $A$ is a
monoid in $\Hmod$, and $L$ is an $A$\h module object (because $\Hmod$ is an
abelian category, we may refer to the monoid $A$ as a ring).  The morphisms
in this global category consist of pairs,
\[(\alpha, f):(L,A) \longrightarrow (M,B)\] 
where $f:A \rightarrow B$ is a ring map and $\alpha:L \rightarrow f^*M$ is a
map of $A$\h modules.  Note that the ring map, $f$, induces a bijection
between the $A$\h module maps and $B$\h module maps:
\[ \frac{L \rightarrow f^*M}{B\otimes_A L \rightarrow M}. \]

Recall that a functor $\pi:\bE \rightarrow \bB$ is said to be a
\definition{fibration} if for every map $f:A \rightarrow B$ in the base
$\bB$, and every object $N$ in $\bE$ which maps down to $B$, we get a
(cartesian) lift of $f$ to $\bE$, $f':M \rightarrow N$, such that given any
map $k:M' \rightarrow N$ where $M'$ is in the fiber over $A$, there exists a
unique map from $M'$ to $M$ making the triangle commute, and which maps down
to the identity morphism on $A$ in $\bB$.  The composites of these
(cartesian) liftings are also a (cartesian) lifting.
\[\xymatrix{
\bE \ar[dd]^{\pi}  & & M  \ar@{-->}[r]^{f'}   & N \\
                   & & M' \ar[ur]_{k} \ar@{-->}[u]        \\ 
\bB                & & A  \ar[r]^{f}    & B
}\]
A \definition{cofibration} is defined dually.  For more information see \cite{borceux}.

\begin{lem}
The functor $\pi:\bE \rightarrow \bB$, from the global category of rings and
modules to the category of monoids (given by $\pi(M,A)=A$) is both a
fibration and a cofibration.  This is often called a
\definition{bifibration}.
\end{lem}

\begin{pf}
We see that the category of rings and modules is a fibration because given
any map of rings, $f:A \rightarrow B$, and any $B$\h module, $M$, the
$A$\h module, $f^*M$, gives us the lift of $f$ to $\bE$:
\[ f^*: (f^*M, A) \longrightarrow (M,B). \]
All maps to $(M,B)$ which project down to $f:A \rightarrow B$ will be of the
form $(\alpha, f):(N,A)\rightarrow (M,B)$, where $\alpha:N \rightarrow f^*M$
is just an $A$\h module map, and hence $\alpha$ just gets mapped to the
identity on $A$.  The proof that this category is also a cofibration goes
through similarly using the adjoint characterization of maps in $\bE$.
\end{pf}

Because we have a bifibration, we may deduce that the global category, $\bE$,
is complete and cocomplete if both the base category is complete and cocomplete,
and each of the fibers is complete and cocomplete.  But both the base and the
fibers are algebras for a monad, and hence are complete and cocomplete.
Note also that by construction, we have a forgetful functor from $\bE$ to
$\Hmod$. 

\section{Maps with Singularities}
\subsection{Binary Singular Maps}
Now that we have the setting of this global category of rings and modules, we
are ready to use its complete and cocomplete structure to form a relaxed
multicategory.  We begin by making precise the notion of singularity we will
be using.

\begin{defn}
An \definition{elementary vertex structure} associated to a (cocommutative)
Hopf algebra $H$ is defined to be an $R$\h module, $K$, which has the
structure of a commutative algebra over $H^*$ as well as that of a 2-sided
$H$\h module.  We require the natural map from $H^*$ to $K$ to be an $H$\h
module morphism.  The algebra structure of $K$ is $H$\h linear, and the
antipode defined on $H^*$ extends to a map from $K^{op}$ to $K$.
\end{defn}

This definition is due to Richard Borcherds, and can be found together with a
number of examples in \cite[Definition 3.2]{bor}.  Intuitively we think of
$K$ as some collection of singular maps defined on $H$.  The following
example motivates the treatment which follows.

\begin{exmp}\label{ex:classical}
Take $H = R[\D{0}, \D{1}, \D{2} \ldots]$ to be a module over a commutative
ring, $R$, with unit.  We give it the structure of a monoid by defining
multiplication $\D{i}\D{j} = \binom{i+j}{i}\D{i+j}$ and unit $\D{0}$, and we
make it into a Hopf algebra by defining comultiplication $\Delta\D{i} =
\sum_{p+q=i} \D{p} \tens \D{q}$, counit $\epsilon(\D{i}) = \delta_{i,0}$, and
antipode $S(\D{i}) = (-1)^i\D{i}$.  Then $H^* \isom R\ps{x}$ and we can take
$K = R\ps{x}[x^{-1}]$ as our elementary vertex structure.  Then for all $j
\in \Z$ we have $\D{i}x^j = \binom{j}{i}\D{j-i}$ and $S(x^j) = (-1)^j x^j$.
\end{exmp}

For any $H$\h module, $D$, the collection of linear maps $\Hom_R(\H{2}, D)$
has the structure of an $\H{2^*}$ module just as $H^*$ has the structure of
an $H^*$\h module.  Using the dual of the map
\[\xymatrix@R=5pt{
H\tens H \ar[r]^f& H \\
h_1\tens h_2\ar@{|->}[r] &    h_1S(h_2)}\]
we consider $K$ as an $\H{2^*}$\h module, and we can tensor over the dual,
$f^*$, to form the module, $\Hom_R(\H{2}, D) \tens_{f^*} K$.  Throughout the
rest of this paper, all tensors with $K$ will be over $f^*$, so we will leave
them from the subscript.  This is an $\H{2}$ module, and so for $H$\h modules
$A$ and $B$, we can consider the collection of $\H{2}$\h linear maps
\begin{equation}\label{e:binary_singular}
A\tens B \longrightarrow \Hom_R(\H{2}, D) \tens_{f^*} K.
\end{equation}
which we call the \definition{singular maps} from $A\tens B$ to $D$.  

\begin{note}
We are interested in this collection of maps because a linear map from
$A\tens B$ to $D$ can be regarded as an $\H{2}$\h linear map from $A\tens B$
to $\Hom_R(\H{2}, D)$.  Hence we've just ``added singularities'' to linear
maps in a way that depends on the underlying Hopf algebra. 
\end{note}

In order to simplify the notation we will be using to describe these singular
maps, we use labelled trees.  The collection of singular maps from $A\tens B$
to $D$ in equation \eqref{e:binary_singular} will be denoted by either of the
following:
\begin{equation*}\label{binary_tree}
\Multibin(A_H,B_H;D) = \flattwoleafed{A_H}{B_H}{D}{x}{y}
\end{equation*}

\noindent
On the right hand side, the leaves of the tree are labelled by the domain of
our singular maps and the root is labelled by the codomain of the singular
maps.  The singularity can be thought of as appearing at the root, and
depending on the inputs above.  The edges are labelled with dummy variables
which act as placeholders for the two copies of $H$.  We use two distinct
dummy variables in order to be able to distinguish the different actions of
$H$.  The $\H{2}$ linearity of our maps is designated by the subscripts on
the $A$ and $B$, and with example \ref{ex:classical} in mind, we could
emphasize this linearity by saying $\partial_A = \partial_x$ and $\partial_B
= \partial_y$.

\begin{note}
Notice that the symmetry of the tensor product implies that the given tree is
isomorphic to its vertical reflection (in terms of the functions they
represent), and $\Multibin(A,B;D) \isom \Multibin(B,A;D)$.
\end{note}

This collection of maps has a number of associated structures which we will
use for working in the global category of rings and modules.  Firstly, there
is the collection of \definition{nonsingular maps}, $\Hom_{\H{2}}\bigl(A\tens
B, \Hom_R(\H{2}, D)\bigr)$ sitting inside the collection of singular maps.
Secondly, we can remove the requirement that the singular and non-singular
maps be $\H{2}$\h invariant, giving \definition{proto-singular} and
\definition{proto-nonsingular} maps from $A\tens B$ to $D$.  And finally,
these proto-singular and proto-nonsingular maps are modules for the rings
$\Hom_R(\H{2}, R)\tens K$ and $\Hom_R(\H{2}, R)$ respectively.  These are
called the \definition{associated rings}.

All of the collections given are $H$\h modules, and so we could consider the
corresponding $H$\h invariant modules.  We sum this up by pointing out that
the collection of proto-singular maps,
\begin{equation*}
\Hom_{R}\Bigl(A\tens B, \Hom_R(\H{2}, D)\tens K \Bigr),
\end{equation*}
has an action of $H$ at $A$, $B$, and $D$, and what we have been calling the
singular maps are just those proto-singular maps which are invariant under
the action of $H$ at both $A$ and $B$.  Similarly, the $H$\h invariant
singular maps are just those maps invariant under the action of $H$ at $A$,
$B$, and $D$.  Using the notation from above to emphasize this $H$\h action,
the proto-singular maps are denoted $\Multibin(A,B;D)$, where we have removed
the $H$ subscript from the $A$ and $B$ as expected.  

\begin{note}
Because we will be working in a category of rings and modules, we focus on
the proto-singular and proto-nonsingular maps, since they are modules for
their associated rings.  We recover our original singular maps by taking
sufficiently $H$\h invariant subcollections.
\end{note}

We now consider composing two proto-singular maps.  Since our motivation was
provided by ordinary multilinear composition, we would like our treatment to
reduce to ordinary linear composition when $K = H^*$.  It is easy to check
that this means requiring the proto-singular maps to be $H$\h invariant at
the point of composition.  So we compose and element of
$\Multibin(A_1,A_2;{B_1}_H)$ with an element of $\Multibin({B_1}_H,B_2;D)$:
\[\begin{split}
A_1 \tens A_2 \longrightarrow  \Hom_H(\H{2}&,  B_1) \tens K\\
& B_1\tens B_2 \longrightarrow \Hom_R(\H{2}, D) \tens K
\end{split}\]
gives an element of 
\begin{equation}\label{e:composed_binaries}
\Hom_{R}\biggl(A_1\tens A_2, \Hom_{H\tens R}\Bigl(\H{2} \tens
B_2, \Hom_R(H\tens H, D) \tens K \Bigr)\tens K \biggr)
\end{equation}
where the inner subscript, $H\tens R$, emphasizes that the only $H$\h
linearity is between the copy of $\H{2}$ tensored with $B_2$ and the first of
the innermost $H$'s.  The associated ring is
\[\Hom_{H\tens R}\Bigl(\H{2} \tens H, \Hom_R(H\tens H, R) \tens K \Bigr)\tens K.\]
Notice that in equation \eqref{e:composed_binaries}, we have taken special
care to emphasize that one of the singularities depends only on inputs $A_1$
and $A_2$, while the other singularity depends on all inputs.  This
dependence of singularities on inputs will be important for the relaxed
multicategory structure we are constructing.  The collection of proto-singular
functions in equation \eqref{e:composed_binaries} also has a corresponding
collection of nonsingular maps, and together with its ring, defines an
inclusion map in the global category.

The collection of proto-singular maps given in equation
\eqref{e:composed_binaries} will be denoted
\[\Multi_{\epsfig{file=Images/3_left.eps,height=3mm}}(A_1, A_2, B_2;D) = \threeleafedleft{A_1}{A_2}{B_2}{D}{x_1}{x_2}{z_1}{z_2}\]

\noindent
The advantage of this notation is that we can see where we have actions of
$H$.  As before we have it between our three inputs, and their corresponding
copies of $H$ (marked with dummy variables $x_1, x_2$ and $z_2$).  The
requirement that the maps from $A_1\tens A_2$ to $B_1$ be $H$\h invariant
at $B_1$ means that the action of $H$ at $z_1$ passes through to an action on $\H{2}$
at $x_1$ and $x_2$.  We use the following suggestive notation to denote this
linearity: $\partial_{x_1} + \partial_{x_2} = \partial_{z_1}$.  Isomorphic
collections of maps appear if we reflect the tree at the internal nodes.

From the discussion of composition and labelled tree notation, it is clear
that by taking suitably $H$\h invariant proto-singular maps we could compose
another binary proto-singular map at either $A_1$ or $A_2$.  Repeating this
process, we see directly how to build up a tower of these proto-singular maps with
only one internal node at each level.  But we should also be able to compose
at $B_2$.  Composing a proto-singular map, $\Multibin(A_3,A_4, {B_2}_H)$,
say, with a map from
$\Multi_{\epsfig{file=Images/3_left.eps,height=3mm}}(A_1, A_2, {B_2}_H;D)$,
we end up with an element of
\begin{equation}\label{e:4comp_1}  
\Hom \Biggl(A_3 \tens A_4, K \tens \Hom \biggl(A_1 \tens A_2, K \tens
\Hom\Bigl(\H{4}, K \tens \Hom(\H{2}, D)\Bigr)\biggr)\Biggr).
\end{equation}
(We have removed the subscripts denoting $H$\h linearity in order to focus
the discussion on the singularities.  We return to the question of linearity
at the end of this section.)  But this collection contains maps which do not appear as
composites of binary proto-singular functions, because there is no singularity here
which depends only on inputs $A_1$ and $A_2$.  Reversing the order of
composition, the map we end up with is an element of
\begin{equation} \label{e:4comp_2} 
\Hom \Biggl(A_1 \tens A_2, K \tens \Hom \biggl(A_3 \tens A_4, K \tens
\Hom\Bigl(\H{4}, K \tens \Hom(\H{2}, D)\Bigr)\biggr)\Biggr).
\end{equation}
These two collections of maps differ only in their dependency of
singularities on inputs.  In particular, the corresponding collection of
nonsingular functions is the same, and so we have inclusion maps in the
global category from the collection of nonsingular functions to both of
\eqref{e:4comp_1} and \eqref{e:4comp_2}.

\begin{note}
We also see that the associated rings are isomorphic and are given by
\begin{equation}\label{e:4ring}
K ^{\otimes 2} \tens \Hom\Bigl(\H{4}, K \tens \Hom(\H{2}, R)\Bigr), 
\end{equation}
where the first of the outer two copies of $K$ is tensored over $f^*$ with
the first two copies of $H$ in $\H{4}$, and similarly for the second copy.
\end{note}

In order to give an exact description of the collection of composites of
these three proto-singular maps,
\begin{equation*}
\flattwoleafed{A_1}{A_2}{{B_1}_H}{x_1}{x_2} \ \
\flattwoleafed{A_3}{A_4}{{B_2}_H}{y_1}{y_2} \ \ 
\flattwoleafed{{B_1}_H}{{B_2}_H}{D}{z_1}{z_2}
\end{equation*}

\noindent
we first take the pushout in the global category of \eqref{e:4comp_1} and
\eqref{e:4comp_2} along the maps including nonsingular functions, and then we
take the pullback over that pushout.  Explicitly, let $T$ denote the ring in
equation \eqref{e:4ring}.  Then in the the fiber of $T$\h modules, we are
taking the pullback of the pushout of the inclusion of
\[ T \tens_{S} \Hom \biggl(A_1 \tens A_2 \tens A_3 \tens A_4,
\Hom\Bigl(\H{4}, \Hom(\H{2}, D)\Bigr)\biggr)\] 
in \eqref{e:4comp_1} and \eqref{e:4comp_2}, where $S$ is the ring associated
to the nonsingular functions.  We denote this collection by the following
labelled tree:
\begin{equation}\label{t:four_double}
\fourleafdoubleC{A}{D}{x}{y}{z}
\end{equation}

\noindent
Throughout this treatment, we have deliberately postponed the discussion of
any $H$\h linearity.  With this labelled tree, it becomes much simpler to see
the actions of our Hopf algebra.  As usual, we have $H$\h actions between the
inputs and the corresponding copies of $H$.  The $H$\h invariance at $B_1$
and $B_2$ adds a further $\H{2}$\h invariance which we denote $\partial_{x_1}
+ \partial_{x_2} = \partial_{z_1}$ and $\partial_{y_1} + \partial_{y_2} =
\partial_{z_2}$ as above.

\subsection{Multimaps Parameterized by Binary Trees}\label{ss:multimaps_binary}
Using the fact that we can represent our proto-singular maps as labelled trees, we
may describe all possible composites of the binary proto-singular maps by
defining proto-singular maps for each binary tree.  The general situation
will be similar to the situation for the tree in \eqref{t:four_double}.  What
will follow will be an algorithmic procedure for describing the
proto-singular maps associated to any binary labelled tree.

For an arbitrary binary tree, we consider its collection of internal
vertices.  We will assume that these include the root, but they do not
include the leaves.  Considering them as a set, this set inherits a partial
order from the tree, where the root is the least element.  We know that any
partial order can be extended to at least one total ordering, possibly many.

Up to this point, our trees have been labeled with $H$\h modules at their
leaves and root.  It will be useful for the explanation which follows to
assume that every internal node is also labelled.  For any internal node,
$q$, connected to $n$ \definition{incoming} nodes (i.e., non-empty nodes
whose height is equal to the height of $q$ plus one and connected to $q$ by a
single edge), we will label $q$ by $\H{n}$.  We can also associate to $q$
the tensor product of the labels of the incoming nodes, and denote it $X_q$.
Thus the following labelled tree has two internal nodes,
\[\threeleafedleft{A_1}{A_2}{A_3}{C}{x_1}{x_2}{z_1}{z_2}\]

\noindent
and we have $X_{\text{root}}=\H{2}\tens A_3$, and
$X_{\text{internal}}=A_1\tens A_2$.

We put our tree in \definition{augmented} form by adding an additional edge
and vertex to the root of our tree.  The root of this new tree is labelled by
the output of the original tree, and the original root is labelled, just as
any internal node, by its inputs.  For binary trees, the original root is
labelled $\H{2}$ because it had two incoming nodes.  We denote the new vertex
$\bottom$, and we therefore have $X_\bottom = \H{2}$.

\begin{defn}\label{d:bin_sing_maps}
Let $p$ be a binary tree with $n$ leaves, and let $t$ denote a total ordering,
$\bottom < \text{root} < p_1 < \cdots < p_l$, of the internal vertices of augmented
$p$, compatible with the the partial ordering inherited from the tree
structure of $p$.  We define an operator on $\H{2}$\h modules:
\begin{equation}
\Sing_{p_i} = \Hom(X_{p_i}, K\tens\cdot).
\end{equation}
Iterating this operator and imposing $H$\h linearity at all internal nodes
(but not the output node), we have
\begin{equation}
\mathrm{Ord}_t(A_1, \ldots, A_n; C) = \Sing_{p_l}\cdots \Sing_{p_1}
\Sing_{\text{root}} \Hom(X_\bottom,C).
\end{equation}
\end{defn}

Notice that for all total orderings, $t$, the collections of nonsingular
functions associated to the $\mathrm{Ord}_t$ are isomorphic, so we are led to
the following definition of $\Multi_p$.

\begin{defn}
$\Multi_p(A_1, \ldots, A_n; C)$ is defined to be the pullback of the pushout
of each $\mathrm{Ord}_t$ for all possible total orderings, $t$, of the
internal vertices of augmented $p$, over the nonsingular functions in the
global category of rings and modules.  
\end{defn}

\begin{note}
As above, the symmetry of
the tensor product implies that vertical reflection of labelled (sub) trees
induces isomorphisms of multimaps.
\end{note}

\begin{exmp}
If $p$ is a tree with only one binary vertex at each level, then there is only
one total ordering, $t$, of internal vertices of the tree, and so 
\[\Multi_p(A_1,\ldots, A_n; C) = \mathrm{Ord}_t(A_1, \ldots, A_n; C).\]
\end{exmp}

\begin{exmp}
When $p = \epsfig{file=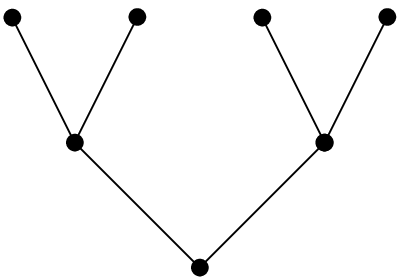,height=3mm}$, there are exactly two
total orderings of internal vertices of this tree, and the corresponding
$\mathrm{Ord}_t(A_1, \ldots, A_4, C)$ functions are given by
\[\begin{split}
\mathrm{Ord}_{t_1} &= \Hom\Bigl(A_1\tens A_2,K\tens
\Hom\bigl(A_3\tens A_4,K\tens \Hom(\H{4},
K\tens\Hom(\H{2},C))\bigr)\Bigr) \\
\mathrm{Ord}_{t_2} &= \Hom\Bigl(A_3\tens A_4,K\tens
\Hom\bigl(A_1\tens A_2, K\tens \Hom(\H{4},
K\tens\Hom(\H{2},C))\bigr)\Bigr),
\end{split}\]
as in equations \eqref{e:4comp_1} and \eqref{e:4comp_2}.  The pullback is
exactly the one described explicitly in the previous section.
\end{exmp}

We finish this section with a proof that composition holds for binary trees.
In order to give this proof, we first need a lemma about evaluation in
symmetric monoidal categories.

\begin{lem}\label{lem:assoc_comp}
In any symmetric monoidal category, $\cC$, the following diagram commutes:
\[\xymatrix{ \Hom(A_1, B_1\tens C_1)\tens\Hom(A_2,B_2\tens C_2)\ar[r]\ar[d]& \Hom(A_1,
C_1 \tens \Hom(A_2,B_1\tens B_2\tens C_2))\ar[d]  \\
\Hom(A_2, C_2 \tens \Hom(A_1,B_1\tens C_1\tens B_2))\ar[r]&
\Hom(A_1\tens A_2, B_1\tens C_1\tens B_2\tens C_2)}\]
\end{lem}

\begin{pf}
The proof follows immediately from the fact that the evaluation of 
\[\Hom(A_1, B_1\tens C_1)\tens\Hom(A_2,B_2\tens C_2)\] 
on $A_1\tens A_2$ gives the same result when carried out by either first
evaluating $A_1$, or by first evaluating $A_2$ or by evaluating both
together.
\end{pf}

\begin{prop}\label{prop:bin_comp}
There exists an associative composition map 
\[\begin{split}
\Multi_q({B_1}_H, \ldots, B_m;& C)  \tens \Multi_{p}(A_1, \ldots, A_n; {B_1}_H)
\longrightarrow \\
&\Multi_{q\circ p}(A_1, \ldots, A_n, B_1, \ldots, B_m; C)
\end{split}\]
for all $H$\h modules $A_i, B_j, C$ and all binary trees $p,q$.
\end{prop}

\begin{pf}
Given two proto-singular maps, $f \in \Multi_{p}(A_1, \ldots, A_n; {B_1}_H)$
and $g\in \Multi_q({B_1}_H, \ldots, B_m; C)$ we prove that they compose to
give an element of $\Multi_{q\circ p}(A_1, \ldots, A_n, B_1, \ldots, B_m;
C)$.  We know that $\Multi_{q\circ p}$ is defined to be a pullback over
objects of the form $\mathrm{Ord}_t(A_1, \ldots, A_n, B_2, \ldots, B_m; C)$,
so we first show that $f$ and $g$ compose to give an element of any such
$\mathrm{Ord}_t$.  This follows from the fact that the linear ordering, $t$,
of the internal nodes of $q\circ p$, induces linear orderings on the internal
nodes of $p$ and $q$.  Denote these linear orderings $t_p$ and $t_q$.
Regarding $f$ and $g$ as elements of $\mathrm{Ord}_{t_p}$ and
$\mathrm{Ord}_{t_q}$, we know that $f\circ g$ is an element of
$\mathrm{Ord}_t$ because of the $H$\h invariance at $B_1$.  By lemma
\ref{lem:assoc_comp}, we know that each image of $f\circ g$ in
$\mathrm{Ord}_t$ gets mapped to the same point in the pushout, and hence they
compose to give an element of $\Multi_{q\circ p}(A_1, \ldots, A_n, B_1,
\ldots, B_m; C)$ as desired.  Associativity of this map is clear.
\end{pf}

\subsection{Non-branching Trees}\label{ss:non_branch_trees}
In order to build up our relaxed multicategory, we want to extend our
definition of multimaps to more general trees.  We begin by expand
definition \ref{d:bin_sing_maps} to allow for non-branching trees:  

\begin{defn} \label{d:1_sing_maps}
If the tree $p$ in definition \ref{d:bin_sing_maps} is allowed to also have
non-branching subtrees, then $\Multi_p(A_1, \ldots, A_n;C)$ is defined
exactly as in that definition except that when an internal vertex, $p_i$ has
only one incoming edge, we define an operator to act on $H$\h modules,
\[\Sing_{p_i} = \Hom_R(X_{p_i},\cdot),\]
where $X_{p_i}$ is the label of the incoming node as in section
\ref{ss:multimaps_binary}.  
\end{defn}

The first and most obvious such tree consists of just a root, $\bullet$.
Composing with such a tree leaves the tree unchanged.  So we hope that this
definition defines $\Multi_\bullet$ so that it composes with a proto-singular
map of type $p$ (for some tree, $p$) to give a proto-singular map of type $p$.
Indeed, labelling the tree $\bullet$ with input and output $H$\h modules, the
definition gives
\[\Multi_\bullet(A;B) = \Hom_R(A,B).\]

Thus we have included all maps from the underlying category.  If we apply
definition \ref{d:1_sing_maps} to the augmented version of the tree,
$\bullet$, we find
\[\Multi_{\ \epsfig{file=Images/1_flat.eps,height=3mm}}(A;B) = \Hom_H(A,\Hom_R(H,B)),\]
and this process can be repeated for any non-branching tree.  Note that
because of the internal $H$\h invariance, we have $\Multi_{\
\epsfig{file=Images/1_flat.eps,height=3mm}}(A;B) \isom \Multi_\bullet(A;B)$.

\begin{exmp}\label{ex:binary_with_tail}
Consider the proto-singular maps associated to the tree,
\epsfig{file=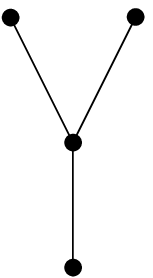,height=3mm}.  From the definition we have
\begin{eqnarray*}
\Multi_{\epsfig{file=Images/2_tail.eps,height=3mm}}(A_1,A_2;B) 
&=& \Hom_{R}\Bigl(A_1\tens A_2, K\tens\Hom_H\bigl(\H{2},\Hom(H,B)\bigr)\Bigr)\\
&\isom& \Hom_{R}\Bigl(A_1\tens A_2, K\tens\Hom\bigl(\H{2},B\bigr)\Bigr)\\
&=& \Multibin(A_1,A_2;B).
\end{eqnarray*}
When $H$ is the classical vertex group this says that there is a bijection
between the collection of proto-singular multimaps associated to the
following trees:
\begin{equation}\label{e:binary_with_tail}
\twoleafedtail{A_1}{A_2}{B}{x}{y}{z} \isom \flattwoleafed{A_1}{A_2}{B}{x+z}{y+z}
\end{equation}

\noindent
This isomorphism follows immediately from the $H$\h invariance at the
internal node, where it provides the relation $\partial_x +\partial_y =
\partial_z$.  
\end{exmp}

Now that we have a description of unary proto-singular maps, it is natural to
consider the nullary type multimaps.  Applying the algorithmic definition for
associating proto-singular maps to trees, we first augment the empty tree, giving
\epsfig{file=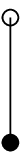,height=4mm}.  As with the tree $\bullet$, the
only internal node is $\bottom$, and we take $X_\bottom = R$ since we do not
consider the empty node as a leaf.  This gives:

\begin{defn}\label{d:0_sing_maps}
For any $H$\h module, $A$, the proto-singular multimaps parameterized by the
empty tree are given by:
\[ \Multi_{\circ}(R;A) = \Hom_R(R,A) \isom A.\]
\end{defn}

What happens when we compose an element of this collection with a binary
proto-singular multimap?

\begin{lem}
We have a composition map
\begin{equation}\label{e:null_comp}
\Multi_{H,\circ}(R;{A_1}_H) \tens \Multibin({A_1}_H,A_2;B) \longrightarrow
\Multi_{\ \epsfig{file=Images/1_flat.eps,height=2mm}}(A_2;B).
\end{equation}
In fact, we have such a composition map for composition with any $\Multi_p$.
\end{lem}

\begin{pf}
Given any binary proto-singular map $f \in \Multibin({A_1}_H,A_2;B)$ and an element $a
\in H\h \text{inv}(A_1) = \Multi_{\circ}(R;{A_1}_H)$, we have a map 
\[f(a\tens \cdot):A_2 \longrightarrow K \tens \Hom_R(\H{2}, B)\]
such that $\epsilon(h)f(a\tens \cdot) = f(ha\tens \cdot)$.  But $f$ is
$H$\h linear at $A_1$, so the map $f$ must factor through 
\[f(a\tens \cdot):A_2 \longrightarrow \Hom_R(H, B),\]
and so we have an element of $\Multi_{\ \epsfig{file=Images/1_flat.eps,height=2mm}}(A_2;B)$.
\end{pf}

\subsection{Ternary Maps and Beyond}
We would now like to define proto-singular maps associated to more general trees.
With the goal of forming a relaxed multicategory, we would like to give a
definition for $\Multi_{\epsfig{file=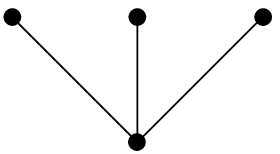,height=3mm}}(A_1, A_2,
A_3;B)$ together with maps to each of the multimaps
$\Multi_{\epsfig{file=Images/3_left.eps,height=3mm}}(A_1, A_2, A_3;B)$,
$\Multi_{\epsfig{file=Images/3_left.eps,height=3mm}}(A_2, A_3, A_1;B)$ and
$\Multi_{\epsfig{file=Images/3_left.eps,height=3mm}}(A_3, A_1, A_2;B)$.  We
know each of these three modules has the same associated collection of
nonsingular functions, together with inclusion maps into each of them, so we
can pullback the pushout of these three objects over the nonsingular
functions in the global category.  This gives us an object which we take as
$\Multi_{\epsfig{file=Images/3_flat.eps,height=3mm}}(A_1, A_2, A_3;B)$,
together with the desired maps.  More generally, we have the following
definition:

\begin{defn}\label{d:n_flat_sing_maps}
For $H$\h modules, $A_1, \ldots, A_n, B$, the collection of proto-singular maps
associated to the flat tree with $n$ leaves,
\epsfig{file=Images/generic_flat.eps,height=3mm}, is denoted
$\Multi_{\epsfig{file=Images/generic_flat.eps,height=3mm}}(A_1, \ldots, A_n;
B)$, and is defined by first taking the pushout of
\[\mathrm{Ord}_t(A_{\sigma(1)}, \ldots, A_{\sigma(n)};B)\]
for all permutations, $\sigma$, and for each total ordering, $t$, of the
internal vertices of binary trees, $p$. with $n$ leaves,
height $n-1$ or less, and no non-branching nodes, over the corresponding
collection of nonsingular functions
\[\Hom\biggl(A_1 \tens \cdots \tens A_n, \Hom\Bigl(\H{2}, \cdots\Hom(\H{2}, B)\cdots\Bigr)\biggr),\]
and then pulling back over this pushout.  
\end{defn}

The idea of this definition is that we take all possible (non-trivial)
collections of proto-singular maps associated to trees which refine to the
flat $n$\h leafed tree, we take the pushout in order to patch the
singularities along the non-singular maps, and we pullback to give an
``intersection'' of the the modules of singularities.  In fact, this
definition suggests a general definition for proto-singular maps associated to
arbitrary trees which generalizes definitions \ref{d:bin_sing_maps},
\ref{d:1_sing_maps} and \ref{d:n_flat_sing_maps}.

\begin{note}
Here we again have isomorphisms of multimaps induced by permutation of input
labels.  
\end{note}

\begin{defn}
For $H$\h modules, $A_1, \ldots, A_n, B$, the collection of proto-singular maps
associated to an arbitrary tree, $q$ with $n$ leaves, $\Multi_q(A_1, \ldots,
A_n; B)$ is as in definition \ref{d:n_flat_sing_maps} except that we pushout
and pullback only those $\mathrm{Ord}$ which can be mapped to from
$\Multi_q(A_1, \ldots, A_n; B)$ by the refinement maps of the relaxed
multicategory.
\end{defn}

\begin{note}
We have symmetries of these general multimaps induced by the symmetries of
the subtrees.
\end{note}

\section{Relaxed Multicategory Structure}
Now that we know how to define $\Multi_p(A_1, \ldots, A_n; B)$ for objects
$A_i, B$, and each $n$\h leafed tree $p$, we have a relaxed multicategory by
taking the fully $H$\h invariant elements of each collection.  In other
words, the multimaps are $\Multi_p({A_1}_H, \ldots, {A_n}_H; {B}_H)$.  

In order to prove that we have defined a relaxed multicategory, we need to
check that we have satisfied the axioms given in definition
\ref{defn_rel_multicat}.  We have satisfied the identity and naturality
axioms by drawing on the underlying categorical structure of $\Hmod$, so we 
only need to show that composition and refinement axioms are satisfied.  We
sketch the proof here.

\begin{thm}
There exists an associative composition map 
\[\begin{split}
\Multi_q({B_1}_H, \ldots, {B_m}_H;& {C}_H) \tens \Multi_{p}({A_1}_H, \ldots, {A_n}_H; {B_1}_H) 
\longrightarrow \\
&\Multi_{q\circ p}({A_1}_H, \ldots, {A_n}_H, {B_1}_H, \ldots, {B_m}_H; {C}_H)
\end{split}\]
for all $H$\h modules $A_i, B_j, C$ and all trees $p,q$.
\end{thm}

\begin{note}
Keep in mind that we are composing the trees $p$ and $q$ and not the
augmented trees.  We only use augmented trees for the purpose of describing
their associated proto-singular multimaps.
\end{note}

\begin{pf}
We defined $\Multi_{q\circ p}(A_1,\ldots, A_n, B_2 \ldots, B_m; C)$ to be the
pullback of all collections of multimaps associated to trees which refine to
$q\circ p$.  So choosing an arbitrary such tree we have refinements of both
$p$ and $q$ which map to the corresponding subtrees.  Thus we are left with
showing that binary trees compose appropriately, which we saw in proposition
\ref{prop:bin_comp}.  It takes a little work to see that each of these
composites is mapped to the same element of the pushout, but is a
straightforward calculation.  So we see that the proto-singular maps do compose to
give an element of $\Multi_{p\circ q}({A_1}_H,\ldots, {A_n}_H, {B_2}_H \ldots, {B_m}_H; {C}_H)$.
\end{pf}

From the construction of $\Multi_p$, we already have nearly all our
refinement maps.  The only refinement maps we excluded were those which
mapped to trees with non-branching internal vertices.  By suitable composition
with the following refinement map, we have all the required refinement maps.

\begin{defn}
A \definition{refinement for a singularity} is the map,
\[K \longrightarrow \Hom_H(H,K)\]
which takes any $k \in K$ to the map $f \in \Hom_H(H,K)$ defined by $f(g) = g
\acts k$ for any $g \in H$.  For any other $H$\h module, $A$, and $H$\h
invariant map $\alpha:A \rightarrow H$, we define a refinement for $K$ by
composition:
\[K \longrightarrow \Hom_H(H,K) \xrightarrow{\alpha} \Hom_H(A,K). \]
\end{defn}

\section{Algebra in the Relaxed Multicategory}

\begin{defn}
An \bemph{(associative) algebra} in a relaxed multicategory, $\cB$, consists
of an object $B \in \cB$ and a collection of maps 
\[\{f_{p}\} = \set{f_{p} \in \Multi_{H,p}(B, \ldots, B, B)}{p \in \T{n}, n
\in \N}.\] 
These maps must
satisfy the following axioms:
\begin{enumerate}
\item\label{d:algebra_composition} \bemph{Composition:} If $q\circ (p_1, \ldots,
p_n)$ is the tree formed by gluing the root of each tree $p_i$ to
the $i$th external edge of an $n$ leafed tree, $q$, ($p_i$ possibly empty), then 
\begin{equation}
f_{q\circ (p_1, \ldots, p_n)} = f_q \circ (f_{p_1}, \ldots, f_{p_n}).
\end{equation}
\item \bemph{Unit:} The map $f_\circ:R \rightarrow B$ (where $\circ$ is the
empty tree) defines a unit for the algebra in the sense that for any $n$
leafed tree, $p$, and any $1 \leq k \leq n$,
\[f_p \circ_k f_\circ = f_{p^\prime}\] 
where $\circ_k$ denotes composition at the $k$th leaf of $p$, and  $p^\prime$
is the $n-1$ leaved tree arrived at by removing the $k$th leaf from $p$.
\item \bemph{Refinement:} If $p,q \in \T{n}$ and $p$ is a refinement of $q$,
then $r_{p,q} (f_p) = f_q$ where $r_{p,q}$ is the refinement map given by the
refinement axiom for a relaxed multicategory.
\end{enumerate}
\end{defn}

This is an algebra in the sense that each map $f_p$ defines an ``$n$\h fold
multiplication'' for elements of $B$.  For all $n \in \N$ we denote the
multimap associated to the flat tree with $n$ leaves by $f_n$.  Since
composition of multimaps in $\cB$ is associative, the associativity of
$(B,\{f_{p}\})$ is a consequence of the composition axiom.  Considering
$\bullet$, the 1 leafed tree with zero edges, then since $f_p \circ_k
f_{\bullet} = f_p$ and $f_{\bullet} \circ f_p = f_p$, we see that
$f_{\bullet} = 1_B$.  The algebra defined by $(B,\{f_{p}\})$ is said to be
\definition{commutative} if the multimaps, $\{f_{p}\}$ are invariant
under an appropriate action of the symmetric group.  This notion makes
sense because our relaxed multicategory is symmetric.

\begin{note}
This definition of an algebra is just a functor from the opposite of the
category of trees to $\cB$ where each object $p \in \mathcal{T}$ is mapped to
an element of $\Multi_{p}$.
\end{note}

Traditional vertex algebras, as found in the literature (e.g., \cite{FHL},
\cite{goddard}, \cite{kac}), arise as exactly algebras for the Hopf algebra
and elementary vertex structure defined in example \ref{ex:classical}, over
the ring $\C$.  The ``locality'' axiom is summed up by the refinement map
$f_{\epsfig{file=Images/3_flat.eps,height=3mm}} \rightarrow
f_{\epsfig{file=Images/3_left.eps,height=3mm}}$, the vacuum is given by
$f_\circ$, and operator product expansions can be deduced from
$f_{\epsfig{file=Images/3_flat.eps,height=3mm}}$.  For more details about the
relation to these axiomatic vertex algebras see \cite{cts_equivalent}.

\ack 
I would like to thank Martin Hyland for his help and guidance.  I would
also like to thank the conference organizers for a very interesting and
beautifully run conference.

\bibliography{full_bib} \bibliographystyle{plain}

\end{document}